\documentclass{svmult}
\usepackage{amssymb}

\ifx\undefined\operatorname\let\operatorname\textrm\fi
\newcommand{\cl}{\operatorname{cl}}

\renewcommand{\Re}{\operatorname{Re\,}}
\renewcommand{\Im}{\operatorname{Im\,}}

\begin{document}

\title{Remarks on Free Entropy Dimension.}

\author{Dimitri Shlyakhtenko
\thanks{Research supported by the National Science Foundation.}}

\institute{Department of Mathematics, UCLA, Los Angeles, CA 90095
\texttt{shlyakht@math.ucla.edu}}

\maketitle

\begin{abstract}
We prove a technical result, showing that the existence of a closable
unbounded dual system in the sense of Voiculescu is equivalent to
the finiteness of free Fisher information. This approach allows one
to give a purely operator-algebraic proof of the computation of the
non-microstates free entropy dimension for generators of groups carried
out in an earlier joint work with I. Mineyev \cite{shlyakht-mineyev:freedim}.
The same technique also works for finite-dimensional algebras. 

We also show that Voiculescu's question of semi-continuity of free
entropy dimension, as stated, admits a counterexample. 
We state a modified version of the question, which avoids the
counterexample, but answering which in the affirmative 
would still imply the non-isomorphism of free group factors.
\end{abstract}

\section*{Introduction.}

Free entropy dimension was introduced by Voiculescu \cite{dvv:entropy2,dvv:entropy3,dvv:entropy5}
both in the context of his microstates and non-microstates free entropy.
We refer the reader to the survey \cite{dvv:entropysurvey} for a
list of properties as well as applications of this quantity in the
theory of von Neumann algebras.

The purpose of this note is to discuss several technical aspects related
to estimates for free entropy dimension.

The first deals with existence of {}``Dual Systems of operators'',
which were considered by Voiculescu \cite{dvv:entropy5} in connection
with the properties of the difference quotient derivation, which is
at the heart of the non-microstates definition of free entropy. We
prove that if one considers dual systems of closed unbounded operators
(as opposed to bounded operators as in \cite{dvv:entropy5}), then
existence of a dual system becomes equivalent to finiteness of free
Fisher information. Using these ideas allows one to give a purely operator-algebraic
proof of the expression for the free entropy dimension of a set of
generators of a group algebra in terms of the $L^{2}$ Betti numbers
of the group \cite{shlyakht-mineyev:freedim}, clarifying the reason
for why the equality holds in the group case. We also point out that
for the same reason one is able to express the non-microstates free
entropy dimension of an $n$-tuple of generators of a finite-dimensional
von Neumann algebra in terms of its $L^{2}$ Betti numbers. In particular,
the microstates and non-microstates free entropy dimension is the
same in this case.

The second aspect deals with the question of semi-continuity of free
entropy dimension, as formulated by Voiculescu in \cite{dvv:entropy2,dvv:entropy3}.
We point out that a counterexample exists to the question of semi-continuity,
as stated. However, the possibility that the free entropy dimension
is independent of the choice of generators of a von Neumann algebra
is not ruled out by the counterexample.

\section{Unbounded Dual Systems and Derivations.}

\subsection{Non-commutative difference quotients and dual systems.}

Let $X_{1},\ldots,X_{n}$ be an $n$-tuple of self-adjoint elements
in a tracial von Neumann algebra $M$. In \cite{dvv:entropy5}, Voiculescu
considered the densely defined derivations $\partial_{j}$ defined
on the polynomial algebra $\mathbb{C}(X_{1},\ldots,X_{n})$ generated
by $X_{1},\ldots,X_{n}$ and with values in $L^{2}(M)\otimes L^{2}(M)\cong HS(L^{2}(M))$,
the space of Hilbert-Schmidt operators on $L^{2}(M)$. If we denote
by $P_{1}:L^{2}(M)\to L^{2}(M)$ the orthogonal projection onto the
trace vector $1$, then the derivations $\partial_{j}$ are determined
by the requirement that $\partial_{j}(X_{i})=\delta_{ij}P_{1}$.

Voiculescu showed that if $\partial_{j}^{*}(P_{1})$ exists, then
$\partial_{j}$ is closable. This is of interest because the existence
of $\partial_{j}^{*}(P_{1})$, $j=1,\ldots,n$ is equivalent to finiteness
of the free Fisher information of $X_{1},\ldots,X_{n}$ \cite{dvv:entropy5}.

Also in \cite{dvv:entropy5}, Voiculescu introduced the notion of
a {}``dual system'' to $X_{1},\ldots,X_{n}$. In his definition,
such a dual system consists of an $n$ tuple of operators $Y_{1},\ldots,Y_{n}$,
so that $[Y_{i},X_{j}]=\delta_{ij}P_{1}$, where Although Voiculescu
required that the operators $Y_{j}$ be anti-self-adjoint, it will
be more convenient to drop this requirement. However, this is not
a big difference, since if $(Y'_{1},\ldots,Y'_{n})$ is another dual
system, then $[Y_{i}-Y'_{i},X_{j}]=0$ for all $i,j$, and so $Y_{i}-Y'_{i}$
belongs to the commutant of $W^{*}(X_{1},\ldots,X_{n})$.

Note that the existence of a dual system is equivalent to the requirement
that the derivations $\partial_{j}:\mathbb{C}(X_{1},\ldots,X_{n})\to HS\subset B(L^{2}(M))$
are inner as derivations into $B(L^{2}(M))$. In particular, Voiculescu
showed that if a dual system exists, then $\partial_{j}:L^{2}(M)\to HS$
are actually closable, and $\partial_{j}^{*}(P_{1})$ is given by
$(Y_{j}-JY_{j}^{*}J)1$. However, the existence of a dual system is
a stronger requirement than the existence of $\partial_{j}^{*}(P_{1})$.

\subsection{Dual systems of unbounded operators.}

More generally, given an $n$-tuple $T=(T_{1},\ldots,T_{n})\in HS^{n}$,
we may consider a derivation $\partial_{T}:\mathbb{C}(X_{1},\ldots,X_{n})\to HS$
determined by $\partial_{T}(X_{j})=T_{j}$ \cite{shlyakht:qdim}.
The particular case of $\partial_{j}$ corresponds to $T=(0,\ldots,P_{1},\ldots,0)$
($P_{1}$ in $j$-th place).

\begin{theorem}
\label{thm:unboubdedDual}Let $T\in HS^{n}$ and assume that $M=W^{*}(X_{1},\ldots,X_{n})$.
The following are equivalent:\\
(a) $\partial_{T}^{*}(P_{1})$ exists;\\
(b) There exists a closable unbounded operator $Y:L^{2}(M)\to L^{2}(M)$,
whose domain includes $\mathbb{C}(X_{1},\ldots,X_{n})$, so that $Y1=0$
and $1$ belongs to the domain of $Y^{*}$, and so that $[Y,X_{j}]=T_{j}$.
\end{theorem}
\begin{proof}
Assume first that (b) holds. Let $\xi=(Y-JY^{*}J)1=JY^{*}1$, which
by assumptions on $Y$ makes sense. Then for any polynomial $Q\in\mathbb{C}(X_{1},\ldots,X_{n})$,\begin{eqnarray*}
\langle\xi,P\rangle & = & \langle(Y-JY^{*}J)1,Q1\rangle\\
 & = & \langle[Y,Q]1,1\rangle\\
 & = & Tr(P_{1}[Y,Q])\\
 & = & \langle P_{1},[Y,Q]\rangle_{HS}\\
 & = & \langle P_{1},\partial_{T}(Q)\rangle_{HS},\end{eqnarray*}
since the derivations $Q\mapsto\partial_{T}(Q)$ and $Q\mapsto[Y,Q]$
have the same values on generators and hence are equal on $\mathbb{C}(X_{1},\ldots,X_{n})$.
But this means that $\xi=\partial_{T}^{*}(P_{1})$.

Assume now that (a) holds. If we assume that $Y1=0$, then the equation
$[Y,X_{j}]=T_{j}$ determines an operator $Y:\mathbb{C}(X_{1},\ldots,X_{n})\to L^{2}(M)$.
Indeed, if $Q$ is a polynomial in $X_{1},\ldots,X_{n}$, then we
have\[
Y(Q\cdot1)=[Y,Q]\cdot1-Q(Y\cdot1)=[Y,Q]\cdot1=\partial_{T}(Q)\cdot1.\]
To show that the operator $Y$ that we have thus defined is closable,
it is sufficient to prove that a formal adjoint can be defined on
a dense subset. We define $Y^{*}$ on $Q\in\mathbb{C}(X_{1},\ldots,X_{n})$
by\[
Y^{*}(Q\cdot1)=-(\partial_{T}(Q^{*}))^{*}\cdot1+\partial_{T}^{*}(P_{1}).\]
Hence $Y^{*}\cdot1=\partial_{T}^{*}(P_{1})$ and $Y^{*}$ satisfies
$[Y^{*},Q]=-(\partial_{T}(Q^{*}))^{*}=-[Y,Q^{*}]^{*}.$ 

It remains to check that $\langle YQ\cdot1,R\cdot1\rangle=\langle Q\cdot1,Y^{*}R\cdot1\rangle$,
for all $Q,R\in\mathbb{C}(X_{1},\ldots,X_{n})$. We have:\begin{eqnarray*}
\langle YQ\cdot1,R\cdot1\rangle & = & \langle[Y,Q]\cdot1,R\cdot1\rangle\\
 & = & \langle1,-[Y^{*},Q^{*}]R\cdot1\rangle\\
 & = & \langle1,Q^{*}Y^{*}R\cdot1\rangle-\langle1,Y^{*}Q^{*}R\cdot1\rangle\\
 & = & \langle Q\cdot1,Y^{*}R\cdot1\rangle-\langle1,Y^{*}Q^{*}R\cdot1\rangle.\end{eqnarray*}
Hence it remains to prove that $\langle1,Y^{*}Q^{*}R\cdot1\rangle=0$.
To this end we write\begin{eqnarray*}
\langle1,Y^{*}Q^{*}R\cdot1\rangle & = & \langle1,[Y^{*},Q^{*}R]\cdot1\rangle-\langle1,Q^{*}RY^{*}\cdot1\rangle\\
 & = & \langle[R^{*}Q,Y]\cdot1,1\rangle-\langle R^{*}Q\cdot1,Y^{*}\cdot1\rangle\\
 & = & Tr([R^{*}Q,Y]P_{1})-\langle R^{*}Q\cdot1,\partial_{T}^{*}(P_{1})\rangle\\
 & = & \langle\partial_{T}(R^{*}Q),P_{1}\rangle_{HS}-\langle\partial_{T}(R^{*}Q),P_{1}\rangle_{HS}=0.\end{eqnarray*}
Thus $Y$ is closable.
\end{proof}
\begin{corollary}
Let $M=W^{*}(X_{1},\ldots,X_{n})$. Then $\Phi^{*}(X_{1},\ldots,X_{n})<+\infty$
if and only if there exist unbounded essentially anti-symmetric
operators $Y_{1},\ldots,Y_{n}:L^{2}(M)\to L^{2}(M)$ whose domain
includes $\mathbb{C}(X_{1},\ldots,X_{n})$, and which satisfy $[Y_{j},X_{i}]=\delta_{ji}P_{1}$.
\end{corollary}
\begin{proof}
A slight modification of the first part of the proof of Theorem \ref{thm:unboubdedDual}
gives that if $Y_{1},\ldots,Y_{n}$ exist, then $\partial_{j}^{*}(P_{1})=(Y_{j}-JY_{j}J)1$
and hence $\Phi^{*}(X_{1},\ldots,X_{n})$ (which is by definition
$\sum_{j}\Vert\partial_{j}^{*}(P_{1})\Vert_{2}^{2}$) is finite.

Conversely, if $\Phi^{*}(X_{1},\ldots,X_{n})<+\infty$, then $\partial_{j}^{*}(P_{1})$
exists for all $j$. Hence by Theorem \ref{thm:unboubdedDual}, we
obtain non-self adjoint closable unbounded operators $Y_{1},\ldots,Y_{n}$
so that the domains of $Y_{j}$and $Y_{j}^{*}$ include $\mathbb{C}(X_{1},\ldots,X_{n})$,
and so that $[Y_{j},X_{i}]=\delta_{ji}P_{1}$. Now since $X_{j}=X_{j}^{*}$
we also have $[Y_{j}^{*},X_{i}]=-\delta_{ji}P_{1}^{*}=-\delta_{ji}P_{1}$.
Hence if we set $\tilde{Y}_{j}=\frac{1}{2}(Y_{j}-Y_{j}^{*})$, we
obtain the desired operators.
\end{proof}

\section{Dual systems and $L^{2}$ cohomology.}

Let as before $X_{1},\ldots,X_{n}\in(M,\tau)$ be a family of self-adjoint
elements.

In conjunctions with estimates on free entropy dimension \cite{shlyakht:qdim,shlyakht-mineyev:freedim}
and $L^{2}$ cohomology \cite{connes-shlyakht:l2betti}, it is interesting
to consider the following spaces:\begin{eqnarray*}
H_{0} & = & \cl{\{ T=(T_{1},\ldots,T_{n})\in HS^{n}:} \\ 
& & \qquad \exists Y\in B(L^{2}(M))\ \ [Y,X_{j}]=T_{j}\}.
\end{eqnarray*}
Here $\cl$ refers to closure in the Hilbert-Schmidt topology.
We also consider
\begin{eqnarray*}
H_{1} & = & \textrm{span}\ \cl{\{ T=(T_{1},\ldots,T_{n})\in HS^{n}:\exists
Y=Y^*\ \textrm{unbounded densely defined}}\\
 &  & \qquad \textrm{with $1$ in the domain of }Y,\ [Y,X_{j}]=T_{j},\  j=1,\ldots,n\}.\end{eqnarray*}
Note that in particular, $H_{0}\subset H_{1}$.

One has the following estimates \cite{shlyakht:qdim,connes-shlyakht:l2betti}:\[
\dim_{M\bar{\otimes}M^{o}}H_{0}\leq\delta^{*}(X_{1},\ldots,X_{n})\leq\delta^{\star}(X_{1},\ldots,X_{n})\leq\Delta(X_{1},\ldots,X_{n}).\]
The main result of this section is the following theorem, whose proof 
has similarities with the Sauvageot's theory of quantum Dirichlet forms
\cite{sauvageot:dirichlet}:

\begin{theorem}
\label{thm:H1isH2}$H_{0}=H_{1}$.
\end{theorem}
\begin{proof}
It is sufficient to prove that $H_{0}$ is dense in $H_{1}$. 

Let $T=(T_{1},\ldots,T_{n})\in HS^{n}$ be such that $T_{j}^{*}=T_{j}=[\I A,X_{j}]$,
$j=1,\ldots,n$, with $A=A^*$ a closed unbounded operator
and $1$ in the domain of $A$. 

For each $0<R<\infty$, let now $f_{R}:\mathbb{R}\to\mathbb{R}$ be a 
$C^{\infty}$ function, so that
\begin{enumerate}
\item $f_{R}(x)=x$ for all $-R\leq x\leq R$; 
\item $|f_{R}(x)|\leq R+1$ for all $x$; 
\item the difference quotient $g_{R}(s,t)=\frac{f_{R}(s)-f_{R}(t)}{s-t}$
is uniformly bounded by $2$.
\end{enumerate}
Let $A_{R}=f_{R}(A)$ and let $T_{j}^{(R)}=[\I A_{R},X_{j}]$. Note
that for each $R$, $T^{(R)}=(T_{1}^{(R)},\ldots,T_{n}^{(R)})\in H_{0}$. Hence 
it is sufficient to prove that $T^{(R)}\to T$ in Hilbert-Schmidt norm as
$R\to\infty$.

Let $\mathcal{A}\cong L^{\infty}(\mathbb{R},\mu)$ be the von Neumann
algebra generated by the spectral projections of $A$; hence $A_{R}\in\mathcal{A}$
for all $R$. If we regard $L^{2}(M)$ as a module over $A$, then
$HS=L^{2}(M)\bar{\otimes}L^{2}(M)$ is a bimodule over $\mathcal{A}$,
and hence a module over $\mathcal{A}\bar{\otimes}\mathcal{A}\cong L^{\infty}(\mathbb{R}^{2},\mu\times\mu)$
in such a way that if $s,t$ are coordinates on $\mathbb{R}^{2}$,
and $Q\in HS$, then $sQ=AQ$ and $tQ=QA$ (more precisely, for any
bounded measurable function $f$, $f(s)Q=f(A)Q$ and $f(t)Q=Qf(A)$).
In particular, we can identify, up to multiplicity, $HS$ with $L^{2}(\mathbb{R}^{2},\mu\times\mu)$.

It is not hard to see that then\[
[f(A),X_{j}]=g\cdot[A,X_{j}],\]
where $g$ is the difference quotient $g(s,t)=(f(s)-f(t))/(s-t)$.
Indeed, it is sufficient to verify this equation on vectors in $\mathbb{C}[X_{1},\ldots,X_{n}]$
for $f$ a polynomial in $A$, in which case the result reduces to\[
[A^{n},X_{j}]=\sum_{k=0}^{n-1}A^{k}[A,X_{j}]A^{n-k-1}=\frac{s^{n}-t^{n}}{s-t}\cdot[A,X_{j}].\]
It follows that\[
T_{j}^{(R)}=[A_{R},X_{j}]=[f_{R}(A),X_{j}]=g_{R}(A)\cdot[A,X_{j}]=g_{R}(A)\cdot T_{j}.\]
Now, since $g_{R}(A)$ are bounded and $g_{R}(A)=1$ on the square
$-R\leq s,t\leq R$, it follows that multiplication operators $g_{R}(A)$
converge to $1$ ultra-strongly as $R\to\infty$. Since $HS$ is a
multiple of $L^{2}(\mathbb{R}^{2},\mu\times\mu)$, it follows that
$g_{R}(A)T_{j}\to T_{j}$ in Hilbert-Schmidt norm. Hence $T^{(R)}\to T$
as $R\to\infty$.
\end{proof}
As a corollary, we re-derive the main result of \cite{shlyakht-mineyev:freedim}
(the difference is that we use Theorem \ref{thm:H1isH2} instead of
the more combinatorial argument \cite{shlyakht-mineyev:freedim};
we sketch the proof to emphasize the exact point at which the fact
that we are dealing with a group algebra becomes completely clear):

\begin{corollary}
Let $X_{1},\ldots,X_{n}$ be generators of the group algebra $\mathbb{C}\Gamma$.
Then $\delta^{\star}(X_{1},\ldots,X_{n})=\delta^{*}(X_{1},\ldots,X_{n})=\Delta(X_{1},\ldots,X_{n})=\beta_{1}^{(2)}(\Gamma)-\beta_{0}^{(2)}(\Gamma)+1$,
where $\beta_{j}^{(2)}(\Gamma)$ are the $L^{2}$-Betti numbers of
$\Gamma$.
\end{corollary}
\begin{proof}
(Sketch). We first point out that in the preceding we could have worked
with self-adjoint families $F=(X_{1},\ldots,X_{n})$ rather than self-adjoint
elements (all we ever needed was that $X\in F\Rightarrow X^{*}\in F$). 

By \cite{shlyakht-mineyev:freedim}, we may assume that $X_{j}\in\Gamma\subset\mathbb{C}\Gamma$,
since the dimension of $H_{0}$ depends only on the pair $\mathbb{C}(X_{1},\ldots,X_{n})$
and its trace. 

Recall \cite{connes-shlyakht:l2betti} that $\Delta(X_{1},\ldots,X_{n})=\dim_{M\bar{\otimes}M^{o}}H_{2}$,
where\[
H_{2}=\{(T_{1},\ldots,T_{n})\in HS:\exists Y^{(k)}\in HS\ \textrm{s.t. }[Y^{(k)},X_{j}]\to T_{j}\ \textrm{weakly}\}.\]
By \cite{connes-shlyakht:l2betti}, $\Delta(X_{1},\ldots,X_{n})=\beta_{1}^{(2)}(\Gamma)-\beta_{0}^{(2)}(\Gamma)+1$;
moreover, from the proof we see that in the group case,\[
H_{2}=\cl\{MXM\},\]
where 
\begin{eqnarray*}
X&=&\{(T_{1}X_{1},\ldots,T_{n}X_{n}):T_{j}\in\ell^{2}(\Gamma),\\
& & \qquad  T_{j}\textrm{ is the value of some }\ell^{2}\ \textrm{group cocycle on }X_{j}\},
\end{eqnarray*}
and where we think of $\ell^{2}(\Gamma)\subset HS$ as {}``diagonal
operators'' by sending a sequence $(a_{\gamma})_{\gamma\in\Gamma}\in\ell^{2}(\Gamma)$
to the Hilbert-Schmidt operator $\sum a_{\gamma}P_{\gamma}$, where
$P_{\gamma}$ is the rank $1$ projection onto the subspace spanned
by the delta function supported on $\gamma$.

Let $\mathfrak{F}$ be the space of all functions on $\Gamma$. Since
the group cohomology $H^{1}(\Gamma;\mathfrak{F}(\Gamma))$ is clearly
trivial, it follows that if $c$ is any $\ell^{2}$-cocycle on $\Gamma$,
then $c(X_{j})=f(X_{j})-f(e)$, for some $f\in\mathfrak{F}$. Hence\[
X=\{([f,X_{1}],\ldots,[f,X_{n}]):f\in\mathfrak{F}\}\cap HS^{n}.\]
Since \emph{every element of $\mathfrak{F}$ is automatically an essentially
self-adjoint operator on $\ell^{2}(\Gamma)$}, whose domain includes
$\mathbb{C}\Gamma$ we obtain that\[
MXM\subset H_{1}.\]
In particular, $H_{2}\subset H_{1}$. Hence \[
\dim_{M\bar{\otimes}M^{o}}H_{2}=\beta_{1}^{(2)}(\Gamma)-\beta_{0}^{(2)}(\Gamma)+1\leq\dim_{M\bar{\otimes}M^{o}}H_{1}\leq\dim_{M\bar{\otimes}M_{0}}H_{2},\]
which forces $H_{1}=H_{2}$. Since $H_{0}=H_{1}$, we get that in
the following equation\[
\dim_{M\bar{\otimes}M^{o}}H_{0}\leq\delta^{*}\leq\delta^{\star}\leq\dim_{M\bar{\otimes}M^{o}}H_{2}=\beta_{1}^{(2)}(\Gamma)-\beta_{0}^{(1)}(\Gamma)+1\]
all inequalities are forced to be equalities, which gives the result.
\end{proof}
\begin{corollary}
Let $(M,\tau)$ be a finite-dimensional algebra, and let $X_{1},\ldots,X_{n}$
be any of its self-adjoint generators. Then $\delta^{*}(X_{1},\ldots,X_{n})=\delta^{\star}(X_{1},\ldots,X_{n})=\Delta(X_{1},\ldots,X_{n})=1-\beta_{0}(M,\tau)=\delta_{0}(X_{1},\ldots,X_{n})$.
\end{corollary}
\begin{proof}
As in the proof of the last corollary, we have the inequalities\[
\dim_{M\bar{\otimes}M^{o}}H_{0}\leq\delta^{*}\leq\delta^{\star}\leq\dim_{M\bar{\otimes}M^{o}}H_{2},\]
where\[
H_{2}=\{(T_{1},\ldots,T_{n})\in HS:\exists Y^{(k)}\in HS\ \textrm{s.t. }[Y^{(k)},X_{j}]\to T_{j}\ \textrm{weakly}\}.\]
Since $L^{2}(M)$ is finite-dimensional, there is no difference between
weak and norm convergence; hence $H_{2}$ is in the (norm) closure
of $\{(T_{1},\ldots,T_{n}):\exists Y\in HS\ \textrm{s.t. }T_{j}=[Y,X_{j}]\}\subset H_{0}$;
since $H_{0}$ is closed, we get that $H_{0}=H_{1}$ and so all inequalities
become equalities. Moreover,\[
\dim_{M\bar{\otimes}M^{o}}H_{2}=\Delta(X_{1},\ldots,X_{n})=1-\beta_{0}(M,\tau)\]
(see \cite{connes-shlyakht:l2betti}).

Comparing the values of $1-\beta_{0}(M,\tau)$ with the computations
in \cite{jung-freexentropy} gives $\delta_{0}(X_{1},\ldots,X_{n})=\delta^{*}(X_{1},\ldots,X_{n})$.
\end{proof}

\section{Some Remarks on Semi-Continuity of Free Dimension.}

In \cite{dvv:entropy2,dvv:entropy3}, Voiculescu asked the question
of whether the free dimension $\delta$ satisfies the following semi-continuity
property. Let $X_{j}^{(k)},X_{j}\in(M,\tau)$ be self-adjoint variables,
$j=1,\ldots,n$, $k=1,2,\ldots$, and assume that $X_{j}^{(k)}\to X_{j}$
strongly, $\sup_{k}\Vert X_{j}^{(k)}\Vert<\infty$. Then is it true
that\[
\liminf_{k}\delta(X_{1}^{(k)},\ldots,X_{n}^{(k)})\geq\delta(X_{1},\ldots,X_{n})?\]
As shown in \cite{dvv:entropy2,dvv:entropy3}, a positive answer to
this question (or a number of related questions, where $\delta$ is
replaced by some modification, such as $\delta_{0}$, $\delta^{*}$,
etc.) implies non-isomorphism of free group factors. In the case of
$\delta_{0}$. a positive answer would imply that the value of $\delta_{0}$
is independent of the choice of generators of a von Neumann algebra.

Although this question is very natural from the geometric standpoint,
we give a counterexample, which shows that some additional assumptions
on the sequence $X_{j}^{(k)}$ are necessary. Fortunately, the kinds
of properties of $\delta$ that would be required to prove the non-isomorphism
of free group factors are not ruled out by this counterexample (see
Question \ref{pro:semicontinuity}).

We first need a lemma.

\begin{lemma}
\label{lem:genFn}Let $X_{1},\ldots,X_{n}$ be any generators of the
group algebra of the free group $\mathbb{F}_{k}$. Then $\delta_{0}(X_{1},\ldots,X_{n})=\delta(X_{1},\ldots,X_{n})=\delta^{*}(X_{1},\ldots,X_{n})=\delta^{\star}(X_{1},\ldots,X_{n})=k$.
\end{lemma}
\begin{proof}
Note that by \cite{guionnet-biane-capitaine:largedeviations} we always
have\[
\delta_{0}(X_{1},\ldots,X_{n})\leq\delta(X_{1},\ldots,X_{n})\leq\delta^{*}(X_{1},\ldots,X_{n})\leq\delta^{\star}(X_{1},\ldots,X_{n});\]
furthermore, by \cite{shlyakht-mineyev:freedim}, $\delta^{\star}(X_{1},\ldots,X_{n})=k$.
Since $\delta_{0}$ is an algebraic invariant \cite{dvv:improvedrandom},
$\delta_{0}(X_{1},\ldots,X_{n})=\delta_{0}(U_{1},\ldots,U_{k})$,
where $U_{1},\ldots,U_{k}$ are the free group generators. Then by
\cite{dvv:entropy3}, $\delta_{0}(U_{1},\ldots,U_{k})=k$. This forces
equalities throughout.
\end{proof}
\begin{example}
Let $u,v$ be two free generators of $\mathbb{F}_{2}$, and consider
the map $\phi:\mathbb{F}_{2}\to\mathbb{Z}/2\mathbb{Z}=\{0,1\}$ given
by $\phi(u)=\phi(v)=1$. The kernel of this map is a subgroup $\Gamma$
of $\mathbb{F}_{2}$, which is isomorphic to $\mathbb{F}_{3}$, having
as free generators, e.g. $u^{2}$, $v^{2}$ and $uv$. Let $X_{1}^{(k)}=\Re u^{2}$,
$X_{2}^{(k)}=\Im u^{2}$, $Y_{1}^{(k)}=\Re v^{2}$, $Y_{2}^{(k)}=\Im v^{2}$,
$Z_{1}^{(k)}=\Re uv$, $Z_{2}^{(k)}=\Im uv$, $W_{1}^{(k)}=\frac{1}{k}\Re u$,
$W_{2}^{(k)}=\frac{1}{k}\Im v$. 

Thus if $X_{1}=\Re u^{2}$, $X_{2}=\Im u^{2}$, $Y_{1}=\Re v^{2}$,
$Y_{2}=\Im v^{2}$, $Z_{1}=\Re uv$, $Z_{2}=\Im uv$, $W_{1}=0$,
$W_{2}=0$, then $X_{j}^{(k)}\to X_{j}$, $Y_{j}^{(k)}\to Y_{j}$,
$Z_{j}^{(k)}\to Z_{j}$ and $W_{j}^{(k)}\to W_{j}$ (in norm, hence
strongly).

Note finally that for $k$ finite, $X_{1}^{(k)},X_{2}^{(k)},Y_{1}^{(k)},Y_{2}^{(k)},Z_{1}^{(k)},Z_{2}^{(k)},W_{1}^{(k_{,})},W_{2}^{(k)}$
generate the same algebra as $u^{2},v^{2},uv,\frac{1}{k}v$, which
is the same as the algebra generated by $u$ and $v$, i.e., the entire
group algebra of $\mathbb{F}_{2}$. Hence $\delta(X_{1}^{(k)},X_{2}^{(k)},Y_{1}^{(k)},Y_{2}^{(k)},Z_{1}^{(k)},Z_{2}^{(k)},W_{1}^{(k_{,})},W_{2}^{(k)})=2$
by Lemma \ref{lem:genFn}. Hence\[
\liminf_{k}\delta(X_{1}^{(k)},X_{2}^{(k)},Y_{1}^{(k)},Y_{2}^{(k)},Z_{1}^{(k)},Z_{2}^{(k)},W_{1}^{(k_{,})},W_{2}^{(k)})=2\]
 On the other hand, $X_{1},X_{2},Y_{1},Y_{2},Z_{1,}Z_{2},W_{1},W_{2}$
generate the same algebra as $u^{2},v^{2},uv,0$, i.e., the group
algebra of $\Gamma\cong\mathbb{F}_{3}.$ Hence\[
\delta(X_{1},X_{2},Y_{1},Y_{2},Z_{1},Z_{2},W_{1},W_{2})=3,\]
which is the desired counterexample.
\end{example}
The same example (in view of Lemma \ref{lem:genFn}) also works for
$\delta_{0}$, $\delta^{*}$ and $\delta^{\star}$.

The following two versions of the question are not ruled out by the
counterexample. If either version were to have a positive answer,
it would still be sufficient to prove non-isomorphism of free group
factors:

\begin{question}
\label{pro:semicontinuity}(a) Let $X_{j}^{(k)},X_{j}\in(M,\tau)$
be self-adjoint variables, $j=1,\ldots,n$, $k=1,2,\ldots$, and assume
that $X_{j}^{(k)}\to X_{j}$ strongly, $\sup_{k}\Vert X_{j}^{(k)}\Vert<\infty$.
Assume that $X_{1},\ldots,X_{n}$ generate $M$ and that for each
$k$, $X_{1}^{(k)},\ldots,X_{n}^{(k)}$ also generate $M$. Then is
it true that\[
\liminf_{k}\delta(X_{1}^{(k)},\ldots,X_{n}^{(k)})\geq\delta(X_{1},\ldots,X_{n})?\]
(b) A weaker form of the question is the following. Let $X_{j}^{(k)},X_{j},Y_{j}\in(M,\tau)$
be self-adjoint variables, $j=1,\ldots,n$, $k=1,2,\ldots$, and assume
that $X_{j}^{(k)}\to X_{j}$ strongly, $\sup_{k}\Vert X_{j}^{(k)}\Vert<\infty$.
Assume that $Y_{1},\ldots,Y_{m}$ generate $M$. Then is it true that\[
\liminf_{k}\delta(X_{1}^{(k)},\ldots,X_{n}^{(k)},Y_{1},\ldots,Y_{m})\geq\delta(X_{1},\ldots,X_{n},Y_{1},\ldots,Y_{m})?\]
\end{question}

We point out that in the case of $\delta_0$,
 these questions are actually equivalent to each other
and to the statement that $\delta_0(Z_1,\ldots,Z_n)$ only depends on the
von Neumann algebra generated by $Z_1,\ldots,Z_n$.

Indeed, it is clear that (a) implies (b).  

On the other hand, if we assume that (b) holds, then we can choose
$X_{j}^{(k)}$ to be polynomials in $Y_{1},\ldots,Y_{m}$, so that
$\delta_{0}(X_{1}^{(k)},\ldots,X_{n}^{(k)},Y_{1},\ldots,Y_{m})=\delta_{0}(Y_{1},\ldots,Y_{m})$
by \cite{dvv:improvedrandom}. Hence $\delta_{0}(Y_{1},\ldots,Y_{m})\geq\delta_{0}(X_{1},\ldots,X_{n},Y_{1},\ldots,Y_{m})\geq\delta_{0}(Y_{1},\ldots,Y_{m})$,
where the first inequality is by (b) and the second inequality is
proved in \cite{dvv:entropy3}. Hence if $W^{*}(X_{1},\ldots,X_{n})=W^{*}(Y_{1},\ldots,Y_{m})$,
then one has $\delta_{0}(X_{1},\ldots,X_{n})=\delta_{0}(X_{1},\ldots,X_{n},Y_{1},\ldots,Y_{m})=\delta_{0}(Y_{1},\ldots,Y_{m})$.
Hence (b) implies that $\delta_0$ is the same on any generators of $M$.

Lastrly, if we assume that $\delta_{0}$ is an invariant of the von
Neumann algebra, then (a)
clearly holds, since the value of
 $\delta_{0}(X_{1}^{(k)},\ldots,X_{n}^{(k)})$
is then independent of $k$ and is 
equal to $\delta_{0}(X_{1},\ldots,X_{n},Y_{1},\ldots,Y_{m})$.

\bibliographystyle{plain}

\end{document}